\input amstex\documentstyle{amsppt}  
\pagewidth{12.5cm}\pageheight{19cm}\magnification\magstep1
\topmatter
\title From classes in the Weyl group to strata\endtitle
\author G. Lusztig\endauthor
\address{Department of Mathematics, M.I.T., Cambridge, MA 02139}\endaddress
\thanks{Supported by NSF grant DMS-2153741}\endthanks
\endtopmatter   
\document
\define\mat{\matrix}
\define\endmat{\endmatrix}

\define\sgn{\text{sgn}}

\define\Irr{\text{\rm Irr}}

\define\lb{\linebreak}

\define\op{\oplus}
   
\define\part{\partial}
\define\emp{\emptyset}

\define\m{\mapsto}
\define\do{\dots}

\define\lra{\leftrightarrow}

\define\T{\times}

\define\nl{\newline}
\redefine\i{^{-1}}

\define\ot{\otimes}

\define\tr{\text{\rm tr}}

\redefine\d{\delta}

\define\r{\rho}
\define\s{\sigma}
\redefine\t{\tau}

\redefine\l{\lambda}
\define\z{\zeta}
\define\x{\xi}

\redefine\D{\Delta}

\define\Th{\Theta}

\define\Ph{\Phi}
\define\Ps{\Psi}

\define\kk{\bold k}

\define\CC{\bold C}

\define\NN{\bold N}

\define\QQ{\bold Q}
\define\RR{\bold R}

\define\ZZ{\bold Z}

\define\cb{\Cal B}

\define\ce{\Cal E}

\define\cp{\Cal P}

\define\cu{\Cal U}

\head Introduction\endhead
\subhead 0.1\endsubhead
Let $W$ be a finite Coxeter group with length function $w\m|w|$.
Let $cl(W)$ be the set of
conjugacy classes of $W$ and let $\Irr(W)$ be the set of
(isomorphism clases of) irreducible representations of $W$ over
$\CC$.

In \cite{L15, 4.1} we have defined (assuming that $W$ is a Weyl
group) a map
${}'\Ph:cl(W)@>>>\Irr(W)$ whose image parametrizes the set of
strata in a reductive group $G$ with Weyl group $W$.

\subhead 0.2\endsubhead
The main contribution of this paper is to propose a new (partly
conjectural) definition of the map ${}'\Ph$ (see \S3)
which, unlike the earlier definition, makes
sense even when $W$ is replaced by a noncrystallographic finite
Coxeter group. The fact that such a definition might exist is
suggested by the fact (see \cite{L15, 4.2}) that when $W$ is a
Weyl group, ${}'\Ph$ depends  only on
$W$ as an abstract Weyl group and not on the underlying root
datum.

\subhead 0.3\endsubhead
Some statements in this paper (such as Theorems 3.5 and 4.2)
were obtained using computer calculation based on the CHEVIE
package, see \cite{CH}, \cite{M15}. I thank Gongqin Li for doing the
programming involved in the calculation.

\head 1. Preliminaries\endhead
\subhead 1.1\endsubhead
In this section $W$ is a Weyl group.
Let $\cp$ be the set of prime numbers union with $\{0\}$.

For $r\in\cp$ let $\kk_r$ be an algebraically closed field of
characteristic $r$ and let $G_r$ be a connected
reductive group over $\kk_r$; we assume that all $G_r$ are
of the same type and have Weyl group $W$.

Let $v$ be an indeterminate. Let $H$ be the Iwahori-Hecke
algebra over $\ZZ[v^2]$ attached to $W$. This is the free
$\ZZ[v^2]$-module with basis
$\{T_w;w\in W\}$ and product
defined by $(T_s+1)(T_s-v^2)=0$ if $|s|=1$,
$T_wT_{w'}=T_{ww'}$ if $|ww'|=|w|+|w'|$; thus, $T_1=1$.
Note that $T_w\m(-v^2)^{|w|}T_{w\i}^{-1}$ defines an algebra
automorphism of $H$.

For any $q\in\RR_{>0}$ let $H(q)$ be the associative
$\CC$-algebra obtained from $H$ be specializing $v^2$ to $q$.
It is known that the irreducible $H(q)$-modules are in natural
bijection $E\lra E(q)$ with $\Irr(W)$.

For any $H(q)$-module $M$ we define a new $H(q)$-module $M^!$
as follows. As $\CC$-vector spaces we have $M^!=M$; for $w\in W$,
the action of
$T_w$ on $M^!$ is equal to the action of $(-q)^{|w|}T_{w\i}^{-1}$ on
$M$. For $E\in\Irr(W)$ we have $E(q)^!=(E\ot\sgn)(q)$ where
$\sgn$ is the sign character of $W$.

For $C\in cl(W)$ let $C_{min}$ be the set of elements of minimal
length in $W$.
For $C\in cl(W),E\in\Irr(W)$ there is a unique $A_{C,E}\in\ZZ[v]$
such that $$A_{C,E}|_{v^2=q}=(q-1)^{m(w)}\tr(T_w,E(q))$$
for any $q\in\RR_{>0}$; here $w\in C_{min}$
and $m(w)$ is the multiplicity of the eigenvalue $1$
in the action of $w$ in the reflection representation of $W$.

The fact that $A_{C,E}$ is independent of the choice of $w$ follows
from results in \cite{GP}.

\subhead 1.2\endsubhead
If $r\in\cp-\{0\}$ and $q$ is a power of $r$ let $G(F_q)$
be the fixed point subgroup of a split Frobenius map $G_r@>>>G_r$
relative to the finite subfield $F_q$ with $q$ elements of
$\kk_r$;
let $\cb(F_q)$ be the fixed point set of this Frobenius map on the
variety of Borel subgroups of $G_r$. Let $\CC[\cb(F_q)]$ be the
vector space
of all functions $\cb(F_q)@>>>\CC$; this is naturally an
$H(q)$-module. It can be decomposed as
$\op_{E\in\Irr(W)}V_q(E)\ot E(q)$ where
$V_q(E)$ are naturally (irreducible) $G(F_q)$-modules.
For any $E\in\Irr(W)$ let $R_q(E)$ be the formal
$\QQ$-linear combination
of (irreducible) unipotent representations of $G(F_q)$ which in
\cite{L78,(3.17.1)} (with $\d=1$) is denoted by $R(E)$.
We can write $R_q(E)=\sum_Vc^V_EV$ where $V$ runs over the
unipotent representations of $G(F_q)$ and $c^V_E\in\QQ$.
We set $\D_q(E)=\sum_Vc^V_E\dim(V)$.

For $E'\in\Irr(W)$, we set $c_{E,E'}=c^{V_q(E)}_{E'}$. Then
$c_{E',E}=c_{E,E'}\in\QQ$ is independent of $q$ (see \cite{L84}).
Moreover we have

(a) $c_{E,E'\ot\sgn}=c_{E\ot\sgn,E'}$.
\nl
From the definitions, for $E\in\Irr(W)$ we have
$$\dim V_q(E)=\sum_{E'\in\Irr(W)}c_{E,E'}\D_q(E').\tag b$$
Letting $q\m1$ this implies
$$\dim(E)=\sum_{E'\in\Irr(W)}c_{E,E'}\dim(E').\tag c$$
For $E,E'$ in $\Irr(W)$ we define
$$A'_{E,E'}=c_{E,E'}.$$
From the definitions, for $w\in W$, we see that
$$\tr(T_w,\CC[\cb(F_q)])\text{ is $0$ if $w\ne1$ and is
$\sharp(\cb(F_q))$ if } w=1.\tag d$$

\subhead 1.3\endsubhead
For $E\in\Irr(W)$ there is a unique $\D(E)\in\ZZ[v^2]$
such that $\D(E)|_{v^2=q}=\D_q(E)$ for any $q$ as in 1.2.
There is a unique $h\in\NN[v^2]$ such that

$h|_{v^2=q}=h(q)$
\nl
for any $q$ as in 1.2 where $h(q)=\sharp(\cb(F_q))(q-1)^{m(1)}$.

For $E,E'$ in $\Irr(W)$ we set
$$A''_{E,E'}=\sum_{\ce\in\Irr(W)}
(\sharp W)\i\sum_{y\in W}\tr(y,E)\tr(y,E')\tr(y,\ce\ot\sgn)\D(\ce)/h,$$
an element of $\QQ(v^2)$.

\subhead 1.4\endsubhead
For $r\in\cp$ let $\cu_r$ be the set of unipotent classes in
$G_r$. Let $\cp_*$ be the set of bad primes for $W$.
The map ${}'\Ph$ in 0.1 was obtained
by combining surjective maps $\Ph^r:cl(W)@>>>\cu_r$ ($r\in\cp$)
defined in \cite{L11} and 
and then using Springer's correspondence $\cu_r@>>>\Irr(W)$ (see
\cite{S76} for $r\in\cp-\cp_*$ and \cite{L81} for $r\in\cp$).

\subhead 1.5\endsubhead
Let $\ZZ[v^2][\Irr(W)]$ be the
free $\ZZ[v^2]$-module with basis $\{E;E\in\Irr(W)\}$.
In \cite{LY} it is shown (by refining results of \cite{L11}) that for
$r\in\cp-\cp_*$, the map $\Ph^r$ is a composition
$\Th^r\Ps$ where $\Ps:cl(W)@>>>\ZZ[v^2][\Irr(W)]$ is independent
of $r$ and $\Th^r$ is a surjective map $\Ps(cl(W))@>>>\cu_r$.

Now $\Ps$ can be viewed as a square matrix
$(\Ps_{C,E})_{C\in cl(W),E\in\Irr(W)}$
with entries $\Ps_{C,E}\in\ZZ[v^2]$. (We have
$\Ps(C)=\sum_{E\in\Irr(W)}\Ps_{C,E}E$.)
As shown in \cite{LY, 2.5}, this matrix is 
a product of three square matrices $AA'A''$ where

$A=(A_{C,E})_{C\in cl(W),E\in\Irr(W)}$ is as in 1.1,

$A'=(A'_{E,E'})_{E\in\Irr(W),E'\in\Irr(W)}$ is as in
1.2, and

$A''=(A''_{E,E'})_{E\in\Irr(W),E'\in\Irr(W)}$ is as in 1.3.

\subhead 1.6\endsubhead
For $w\in W$ we have
$(-v^2)^{|w|}T_{w\i}\i=\sum_{y\in W}e_{y,w}T_y$
where $e_{y,w}\in\ZZ[v^2]$.
Note that $e_{y,w}=0$ unless $y\le w$.
For $y,w$ in $W$ let $e_{y,w}(0)\in\ZZ$ be the value of
$e_{y,w}$ at $v=0$.

Assume that $w,s$ in $W$ satisfy $|s|=1,|sw|=|w|+1$. We have
$$\align&
\sum_{y\in W}e_{y,sw}T_y=(-T_s+v^2-1)\sum_{y\in W}e_{y,w}T_y\\&=
-\sum_{y;|sy|>|y|}e_{y,w}T_{sy}-\sum_{y;|sy|<|y|}v^2e_{y,w}T_{sy}
+\sum_{y;|sy|>|y|}(v^2-1)e_{y,w}T_y.\endalign$$
It follows that

(a) $e_{y,sw}=-v^2e_{sy,w}+(v^2-1)e_{y,w} \text{ if }|sy|>|y|$,

(b) $e_{y,sw}=-e_{sy,w} \text{ if }|sy|<|y|$.

In particular,

(c) $e_{y,sw}(0)=-e_{y,w}(0) \text{ if }|sy|>|y|$,

(d) $e_{y,sw}(0)=-e_{sy,w} \text{ if }|sy|<|y|$.

\proclaim{Lemma 1.7}For $y,w$ in $W$ such that $y\le w$ we have
$e_{y,w}(0)=(-1)^{|w|}$.
\endproclaim
We argue by induction on $|w|$. If $w=1$ the result is obvious.
Assume that the result is known for some $w$. We prove it for
$sw$ where $s\in W$ satisfies $|s|=1,|sw|=|w|+1$.
Let $y'\in W$ be such that $y'\le sw$. If $|sy'|>|y'|$ then $y'\le w$
and from 1.6(c) and the induction hypothesis we see that
$e_{y',sw}(0)=(-1)^{|sw|}$ as required. If $|sy'|<|y'|$ then $sy'\le w$
and from 1.6(d) and the induction hypothesis we see that
$e_{y',sw}(0)=(-1)^{|sw|}$ as required.

\proclaim{Lemma 1.8} Let $y,w$ be in $W$.

(i) $v^{2|y|}e_{y,w}$ is of the form
$av^{2|w|}+\text{ lower powers of }v$ for some $a\in\ZZ$;

(ii) if $y=1$ then $a$ in (i) is $1$.
\endproclaim
We prove (i) by induction on $|w|$. If $w=1$ the result is obvious.
Assume that the result is known for some $w$. We prove it for
$sw$ where $s\in W$ satisfies $|s|=1,|sw|=|w|+1$.
Using 1.6(a) we see that
if $|sy|>|y|$ then
$$v^{2|y|}e_{y,sw}=-v^{2|sy|}e_{sy,w}+v^{2|y|}(v^2-1)e_{y,w}.$$
By the induction hypothesis the right hand side
is of the form
$$(-av^{2|w|}+\text{ lower powers of }v)+
(a'v^{2|sw|}+\text{ lower powers of }v)$$
(with $a,a'$ in $\ZZ$) hence the left hand side
is of the form $a'v^{2|sw|}+\text{ lower powers of }v$.
Using 1.6(b) we see that if
$|sy|<|y|$ then  $v^{2|y|}e_{y,sw}=-v^{2|sy|}v^2e_{sy,w}$.
By the induction hypothesis the right hand side
is of the form $-av^{2|sw|}+\text{ lower powers of }v$ (with $a\in\ZZ$)
hence so is the left hand side. This proves (i).

We prove (ii) by induction on $|w|$. If $w=1$ the result is obvious.
Assume that the result is known for some $w$. We prove it for
$sw$ where $s\in W$ satisfies $|s|=1,|sw|=|w|+1$. By 1.6(a) we have

$e_{1,sw}=-v^2e_{s,w}+(v^2-1)e_{1,w}$.
\nl
By the induction hypothesis we have
$$(v^2-1)e_{1,w}=(v^2-1)(v^{2|w|}+\text{ lower powers of } v
=v^{2|sw|}+\text{ lower powers of }v).$$
By (i) we have $v^2e_{s,w}=av^{2|w|}+\text{ lower powers of }v$
where $a\in\ZZ$.
It follows that $e_{1,sw}=v^{2|sw|}+\text{ lower powers of }v$.
This proves (ii).

\proclaim{Lemma 1.9}Let $w\in W$. We have $|w|+m(w)-m(1)=0\mod2$.
\endproclaim
Assume first that $m(w)=0$. Let $\det(-w,\r)$ be the determinant of $-w$ on
the reflection representation $\r$ of $W$. We have
$\det(-w,\r)=(-1)^{|w|-m(1)}$.
It is enough to show that $\det(-w,\r)=1$. It is known that $w,w\i$ are
conjugate in $W$. Hence the eigenvalues of $w:\r@>>>\r$ occur in pairs
$\l\ne\l\i$ and there may be also some eigenvalues $-1$ (there are no
eigenvalues $1$). Thus  the eigenvalues of $w:\r@>>>\r$ occur in pairs
$-\l\ne-\l\i$ and there are also some eigenvalues $1$ (there are no
eigenvalues $-1$). We see that $\det(-w,\r)=1$ as desired.
Next we assume that $m(w)>0$. Then $w$ is conjugate to an element in a
proper standard parabolic subgroup $W'$ of $W$.
Then the desired result for
$w$ follows from the analogous result for $W'$ which can be assumed
known.

\proclaim{Lemma 1.10}Let $w\in W$. We have $|w|+m(w)-m(1)\ge0$.
\endproclaim
Assume first that $m(w)=0$. If $|w|<m(1)$ then $w$ is contained in
a proper standard parabolic subgroup $W'$ of $W$. But then $w$ must
have some eigenvalue $1$ on the reflection representation, contradicting
$m(w)=0$. We see that $|w|+m(w)-m(1)\ge0$.
Next we assume that $m(w)>0$. Then $w$ is conjugate to an element in a
proper standard parabolic subgroup $W'$ of $W$.
Then the desired result for
$w$ follows from the analogous result for $W'$ which can be assumed
known.

\head 2. The polynomials $\Ps_{C,1},\Ps_{C,\sgn},\Ps_{\{1\},E}$
\endhead
\subhead 2.1\endsubhead
In this section $W$ is a Weyl group.
For $\x\in H$ we can write $\x=\sum_{w\in W}\x[w]T_w$
where $\x[w]\in\ZZ[v^2]$. Similarly for $q$ as in 1.2 and
$\x\in H_q$ we can write $\x=\sum_{w\in W}\x[w]T_w$
where $\x[w]\in\CC$.

\proclaim{Proposition 2.2}Let $C\in cl(W)$ and let $w\in C_{min}$.

(i) We have $\Ps_{C,1}=(v^2-1)^{m(w)-m(1)}((-v^2)^{|w|}T_{w\i}\i)[1]$.

(ii) $\Ps_{C,\sgn}$ is $1$ if $C=\{1\}$ and is $0$ if $C\ne\{1\}$.
\endproclaim
Let $q$ be as in 1.2. We prove (i). By 1.5 we have
$$\align&
\Ps_{C,1}|_{v^2=q}=\sum_{E'\in\Irr(W),E''\in\Irr(W),\ce\in\Irr(W)}
A_{C,E'}|_{v^2=q}A'_{E',E''}A''_{E'',1}|_{v^2=q}\\&
=\sum_{E'\in\Irr(W),E''\in\Irr(W),\ce\in\Irr(W)}
(q-1)^{m(w)}\tr(T_w,E'(q))c_{E',E''}\times \\&
\times(\sharp W)\i\sum_{y\in W}\tr(y,E'')\tr(y,\ce\ot\sgn)\D_q(\ce)/h_q.
\endalign$$
Now $(\sharp W)\i\sum_{y\in W}\tr(y,E'')\tr(y,\ce\ot\sgn)$ is $1$ if
$E''=\ce\ot\sgn$ and is $0$ if $E''\ne\ce\ot\sgn$. Hence
$$\Ps_{C,1}|_{v^2=q}=\sum_{E'\in\Irr(W),\ce\in\Irr(W)}
(q-1)^{m(w)}\tr(T_w,E'(q))c_{E',\ce\ot\sgn}\D_q(\ce)/h_q.$$
Using 1.2(a),(b), this becomes
$$\align&\Ps_{C,1}|_{v^2=q}=\sum_{E'\in\Irr(W)}
(q-1)^{m(w)}\tr(T_w,E'(q))\sum_{\ce\in\Irr(W)}c_{E'\ot\sgn,\ce}
\D_q(\ce)/h_q\\&
=\sum_{E'\in\Irr(W)}(q-1)^{m(w)}\tr(T_w,E'(q))
\dim V_q(E'\ot\sgn)/h_q\\&
=\sum_{E_1\in\Irr(W)}(q-1)^{m(w)}\tr(T_w,E_1(q)^!)
\dim V_q(E_1)/h_q\\&
=\sum_{E_1\in\Irr(W)}(q-1)^{m(w)}\tr((-q)^{|w|}T_{w\i}\i,E_1(q))
\dim V_q(E_1)/h_q\\&
=(q-1)^{m(w)}\tr((-q)^{|w|}T_{w\i}\i,\CC[\cb(F_q)])/h_q.\endalign$$
Using 1.2(d) we see that this equals
$$(q-1)^{m(w)}((-q)^{|w|}T_{w\i}^{-1})[1]\sharp(\cb(F_q))/h_q=
(q-1)^{m(w)-m(1)}((-q)^{|w|}T_{w\i}^{-1})[1].$$
Now (i) follows.

The proof of (ii) is similar (but simpler). By 1.5 we have
$$\align&
\Ps_{C,\sgn}|_{v^2=q}=\sum_{E'\in\Irr(W),E''\in\Irr(W),\ce\in\Irr(W)}
A_{C,E'}|_{v^2=q}A'_{E',E''}A''_{E'',\sgn}|_{v^2=q}\\&
=\sum_{E'\in\Irr(W),E''\in\Irr(W),\ce\in\Irr(W)}
(q-1)^{m(w)}\tr(T_w,E'(q))c_{E',E''}\times \\&
\times
(\sharp W)\i\sum_{y\in W}\tr(y,E'')\sgn(y)\tr(y,\ce\ot\sgn)\D_q(\ce)/h_q.
\endalign$$
Now $(\sharp W)\i\sum_{y\in W}\tr(y,E'')\sgn(y)\tr(y,\ce\ot\sgn)$ is
$1$ if $E''=\ce$ and is $0$ if $E''\ne\ce$. Hence
$$\Ps_{C,\sgn}|_{v^2=q}
=\sum_{E'\in\Irr(W),\ce\in\Irr(W)}(q-1)^{m(w)}
\tr(T_w,E'(q))c_{E',\ce}\D_q(\ce)/h_q.$$
Using 1.2(b), this becomes
$$\align&\Ps_{C,1}|_{v^2=q}
=\sum_{E'\in\Irr(W)}(q-1)^{m(w)}\tr(T_w,E'(q))\dim V_q(E')/h_q\\&
=(q-1)^{m(w)}\tr(T_w,\CC[\cb(F_q)])/h_q.\endalign$$
Using 1.2(d) we see that (ii) holds.

\proclaim{Proposition 2.3}Let $E\in\Irr(W)$. We have
$\Ps_{\{1\},E}=\dim(E)$.
\endproclaim
By 1.5 we have
$$\align&\Ps_{1,E}=\sum_{E'\in\Irr(W),E''\in\Irr(W),\ce\in\Irr(W)}
A_{1,E'}A'_{E',E''}A''_{E'',E}\\&
=\sum_{E'\in\Irr(W),E''\in\Irr(W),\ce\in\Irr(W)}
(v^2-1)^{m(1)}\dim(E')c_{E',E''}\times \\&
\times(\sharp W)\i\sum_{y\in W}\tr(y,E'')\tr(y,E)\tr(y,\ce\ot\sgn)\D(\ce)/h.
\endalign$$
Using 1.2(c) we deduce
$$\align&\Ps_{1,E}=\sum_{E''\in\Irr(W),\ce\in\Irr(W)}
(v^2-1)^{m(1)}\dim(E'')\times \\&
\times(\sharp W)\i\sum_{y\in W}\tr(y,E'')\tr(y,E)\tr(y,\ce\ot\sgn)\D(\ce)/h.
\endalign$$
Now $\sum_{E''\in\Irr(W)}\dim(E'') \tr(y,E'')$ is $0$ if $y\ne1$
and is $\sharp W$ if $y=1$. Hence
$$\align&\Ps_{1,E}=\sum_{\ce\in\Irr(W)}
(v^2-1)^{m(1)}\dim(E)\dim(\ce\ot\sgn)\D(\ce)/h\\&
=\sum_{\ce\in\Irr(W)}
(v^2-1)^{m(1)}\dim(E)\dim(\ce)\D(\ce)/h.\endalign$$
From the definitions, for $q$ as in 1.2 we have
$$\sum_{\ce\in\Irr(W)}\dim(\ce)\D_q(\ce)=\sharp(\cb(F_q)).$$
We deduce
$$\Ps_{1,E}|_{v^2=q}=\dim(E)(q-1)^{m(1)}\sharp(\cb(F_q))/h_q=\dim(E).$$
This completes the proof.

\proclaim{Corollary 2.4}Let $C,w$ be as in 2.2.

(i) The leading term of the polynomial $\Ps_{C,1}$ (in $v^2$) is
$v^{2(|w|+m(w)-m(1))}$. (Recall that $|w|+m(w)-m(1)\ge0$, see 1.10.)

(ii) The constant term of the polynomial $\Ps_{C,1}$ (in $v^2$) is $1$.
\endproclaim
By 2.2(i), we have
$\Ps_{C,1}=(v^2-1)^{m(w)-m(1)}e_{1,w} $ (notation of 1.6).
By 1.8(ii) we have $e_{1,w}=v^{2|w|}+\text{ lower powers of }v$. We have
$$(v^2-1)^{m(w)-m(1)}=v^{2(m(w)-m(1)}+\text{ a possibly infinite
sum of lower powers of }v.$$
Hence
$$\Ps_{C,1}=v^{2(|w|+m(w)-m(1))}
+\text{ a possibly infinite sum of lower powers of }v.$$
Since $\Ps_{C,1}$ is a polynomial in $v^2$ the last sum is finite.
This  proves (i).

By 1.7 we have
$e_{1,w}=(-1)^{|w|}+\text{ strictly positive powers of }v$. We have
$$\align&(v^2-1)^{m(w)-m(1)}\\&
=(-1)^{(m(w)-m(1)}+\text{ a possibly infinite
sum of strictly positive powers of }v.\endalign$$
Hence
$$\align&\Ps_{C,1}\\&
=(-1)^{|w|+m(w)-m(1)}
+\text{ a possibly infinite sum of strictly positive
powers of }v.\endalign$$
Since $\Ps(C,1)$ is a polynomial in $v^2$ the last sum is finite.
Thus the polynomial $\Ps(C,1)$ has constant term $(-1)^{|w|+m(w)-m(1)}$.
This equals $1$ by 1.9. This proves (ii).

\head 3. The ``definition'' of $\s:cl(W)@>>>\Irr(W)$\endhead
\subhead 3.1\endsubhead
In this section $W$ is a Weyl group; we give a
(partly conjectural) definition of a map $\s:cl(W)@>>>\Irr(W)$ such
that $\s(C)={}'\Ph(C)$ for any $C\in cl(W)$. This uses the map
$\Ps$ and does not rely on the geometry of $G_r$. 

\subhead 3.2\endsubhead
Let $\Irr_*(W)$ be the set of all $E\in\Irr(W)$ such that
for some $C\in cl(W)$, we have $\Ps_{C,E}\ne0$ and the coefficient
of the highest power of $v$ in $\Ps_{C,E}$ is $<0$.

For $C\in cl(W)$ let $X_C=\{E\in\Irr(W)-\Irr_*(W);\Ps_{C,E}\ne0\}$.
By 2.4, we have $1\in X_C$. In particular $X_C$ is nonempty.
For $E\in\Irr(W)$ let $b_E\in\NN$  be the smallest integer $n$ such
that $E$ appears in the $n$-the symmetric power of the reflection
representation of $W$. Let $X_C^{max}$ be the set of all $E\in X_C$
such that $b_{E'}\le b_E$ for all $E'\in X_C$.

\proclaim{Conjecture 3.3}(a) For any $C\in cl(W)$, $X_C^{max}$
consists of a single element, denoted $\s(C)$. This defines a map
$\s:cl(W)@>>>\Irr(W)$.

(b) We have $\s(C)={}'\Ph(C)$ for all $C\in cl(W)$. The image of
$\s={}'\Ph$ is equal to $\Irr(W)-\Irr_*(W)$.

(c) For any $C\in cl(W)$ we have $\Ps_{C,\s(C)}=v^{|w|+m(w)-m(1)}$
where $w\in C_{min}$.
\endproclaim

\subhead 3.4\endsubhead
If $C=\{1\}$ we have $\Ps_{C,\sgn}=1$ by 2.2(ii) or 2.3. Hence
for this $C$, 3.3(a) holds and we have
$\s(C)=\sgn={}'\Ph(C)$.

In the case where $W$ is of rank $\le2$ the validity of the
conjecture can be deduced from results in \cite{LY}.

\proclaim{Theorem 3.5}Conjecture 3.3 holds when $W$ is of
type $B_3,F_4,E_6,E_7,E_8$.
\endproclaim
For $W$ as in the Theorem, the map $\s={}'\Ph$
is described explicitly in \cite{L15}.

\subhead 3.6\endsubhead
The following errata to \cite{L15} overrides the errata in \cite{L22}.

On p.355, line 6, replace $[A_3+A_2]$ by $[D_4(a_1)+2A_1]$.

On p.356, line 13, replace $[(A_5+A_1)']$ by $[(A_5+A_1)',A_5]$.

\head 4. The noncrystallographic case\endhead
\subhead 4.1\endsubhead
In this section, $W$ denotes an irreducible
noncrystallographic finite Coxeter group.
As mentioned in \cite{LY}, $\Ps$ makes sense even for such $W$.
(The analogue of $A'$ in 1.2 is given for dihedral $W$
in \cite{L94} and for $W$ of type $H_4$ in \cite{M94}.)
But in this case $\Ps$ is only a map $cl(W)@>>>\QQ(v^2)[\Irr(W)]$;
the fact that its image is contained in $\ZZ[v^2][\Irr(W)]$ is a miracle
which is a result of computation.
In any case,
the subset $\Irr_*(W)$ of $\Irr(W)$ and (for $C\in cl(W)$) the
subsets $X_C,X_C^{max}$ of $\Irr(W)-\Irr_*(W)$ can be defined (as in
3.2).
But there is a difference from what was stated in Conjecture 3.3.
In the present case the set $X_C^{max}$ (for $C\in cl(W)$)
may contain more than one element. 

Let $\Irr_{**}(W)$ be the set of all $E\in\Irr(W)-\Irr_*(W)$
such that for some $C\in cl(W)$, we have $E\in X_C^{max}$
and the polynomial $\Ps_{C,E}$ (in $v$) has some $<0$ coefficient.
For $C\in cl(W)$ let
${}'X_C=\{E\in(\Irr(W)-\Irr_*(W))-\Irr_{**}(W);\Ps_{C,E}\ne0\}$.
We have $1\in{}' X_C$. In particular ${}'X_C$ is nonempty.
Let ${}'X_C^{max}$ be the set of all $E\in{}'X_C$
such that $b_{E'}\le b_E$ for all $E'\in{}'X_C$.
(In the setup of 3.3, $\Irr_{**}(W)$ would be empty and
${}'X_C=X_C,{}'X_C^{max}=X_C^{max}$.)

\proclaim{Theorem 4.2}For any $C\in cl(W)$, ${}'X_C^{max}$
consists of a single element, denoted $\s(C)$. This defines a map
$\s:cl(W)@>>>\Irr(W)$.
\endproclaim
Note that the analogue of 3.3(c) is not true in types $H_3,H_4$; see the tables in 4.3, 4.4.

\subhead 4.3\endsubhead
In this subsection we assume that $W$ is of type $H_4$. Assume that
the simple reflections are numbered so that
$(s_1s_2)^5=(s_2s_3)^3=(s_3s_4)^3=1$.
For a sequence $i_1,i_2,\do,i_k$ in $\{1,2,3,4\}$ we write
$(i_1i_2\do i_n)$ for the conjugacy class in $W$ containing
$s_{i_1}s_{i_2}\do s_{i_n}$. We write $c_n$ for an elliptic conjugacy
class in $W$ containing an element of minimal length $n$. (This
is ambiguous only if $n=16$ when there are two such conjugacy
classes, $c_{16},c'_{16}$.
The set $\Irr(W)$ consists of $34$ objects:
$$\align&1_0,4_1,9_2,16_3,25_4,36_5,16_6,9_6,24_6,4_7,24_7,40_8,\\&
48_9,30_{10},30'_{10},18_{10},24_{11},16_{11},24_{12},10_{12},8_{12},\\&
6_{12},16_{13},8_{13},36_{15},25_{16},16_{18},6_{20},16_{21},9_{22},\\&
9_{26},4_{31},4_{37},1_{60}.\endalign$$
(We write $d_b$ for $E\in \Irr(W)$ with $\dim(E)=d$, $b_E=b$. This is
ambiguous only for $d=30,b=10$ when there are two such $E$: $30_{10},
30'_{10}$.)

The table below describes the map $\s$ in the form $?...?...?$
where the first $?$ represent the various $C$ which map to a given
$E$ (represented by the second $?$) and the third $?$ gives the
polynomials $\Ps_{C,E}$ for $C$ in the first $?$ and $E$ in the second
$?$.

$c_4...1_0...v^0$

$c_6...4_1...v^2$

$c_8...9_2 ...v^4$

$c_{10}...16_3...v^6$

$c_{12}...25_4...v^8$

$c_{16},c_{14},(123)...36_5...v^{14}+v^{10},v^{10},v^0$

$c'_{16}...24_6...v^{12}$

$c_{18}...24_7...v^{14}$

$c_{20}...40_8...v^{16}$

$c_{22},(12123)...48_9  ...v^{18},v^2$

$c_{24}...18_{10}...v^{20}$

$c_{26}...24_{11}...v^{22}$

$(124)... 24_{12}...v^0$

$c_{28}...8_{12}...v^{24}$    

$c_{30},(243)...8_{13}...v^{26},v^0$

$(12)...36_{15}...v^0$

$c_{40},c_{38},c_{36},(123)^3,(134)...25_{16}...v^{40}+v^{36}+v^{32},
                        v^{36}+v^{32},v^{32},v^6,v^0$

$(23)...16_{21}...v^0$

$c_{48},(12124)...9_{22}...v^{44},v^2$

$c_{60},(123)^5,(1212),(13)...9_{26}...v^{60}+v^{56}+v^{52},
                        v^{14}+v^{10},v^2,v^0$

$ (1)...4_{37}...v^0$

$(-)...1_{60}....v^0$
                     
  $\Irr_*(W)$ consists of $9$ objects: $4_7,30_{10},16_{11},
  10_{12},6_{12},16_{13},16_{18},6_{20},4_{31}$;
  $\Irr_{**}(W)$ consists of $3$ objects: $9_6,16_6,30'_{10}$.
                        
 The image of $\s$ is equal to $(\Irr(W)-\Irr_*(W))-\Irr_{**}(W)$; it
                        consists of the $22$ objects
 $$\align&
 1_0,4_1,9_2,16_3,25_4,36_5,24_6,24_7,40_8,48_9,18_{10},\\&24_{11},
 24_{12},8_{12},8_{13},36_{15},25_{16},16_{21},9_{22},9_{26},
 4_{37},1_{60}\endalign$$

 \subhead 4.4\endsubhead
In this subsection we assume that $W$ is of type $H_3$. Assume that
the simple reflections are numbered so that
$(s_1s_2)^5=(s_2s_3)^3=1$.
For a sequence $i_1,i_2,\do,i_k$ in $\{1,2,3\}$ we write
$(i_1i_2\do i_n)$ for the conjugacy class in $W$ containing
$s_{i_1}s_{i_2}\do s_{i_n}$. We write $c_n$ for an elliptic conjugacy
class in $W$ containing an element of minimal length $n$.
The set $\Irr(W)$ consists of $10$ objects:
$$1_0,3_1,5_2,3_3,4_3,4_4,5_5,3_6,3_8,1_{15}.$$
(We write $d_b$ for $E\in \Irr(W)$ with $\dim(E)=d$, $b_E=b$.)
The table below describes the map $\s$ in the form $?...?...?$
(with conventions similar to those in 4.3.)

$c_3...1_0...v^0$

$c_5...3_1...v^2$

$(12)...3_3....v^0$

$c_9...4_3...v^6$

$(23)...4_4...v^0$

$c_{15},(1212),(13)...5_5...v^{14}+v^{10},v^2,v^0$

$(1)...3_8...v^0$

$(-)...1_{15}...v^0$.

$\Irr_*(W)$ consists of $1$ object: $3_6$; $\Irr_{**}(W)$ is empty.

 The image of $\s$ is contained strictly in
 $(\Irr(W)-\Irr_*(W))-\Irr_{**}(W)$; it
        consists of the $8$ objects

$1_0,3_1,3_3,4_3,4_4,5_5,3_8,1_{15}$.

\subhead 4.5\endsubhead
The case where $W$ is a finite dihedral group will be
considered in \S5.   

\head 5. Dihedral groups\endhead
\subhead 5.1\endsubhead
In this section $W$ is the dihedral group of order $2p$ where
$p\in\{3,4,5,\do\}$.
It has generators $s,t$ and relations $s^2=t^2=(st)^p=1$.

Here is a list of representatives of the various classes in $cl(W)$ (we set $c=st$):

(a) $c,c^2,\do,c^{(p-1)/2},s,1$ ($p$ odd),

(b) $c,c^2,\do,c^{p/2},s,t,1$ ($p$ even).

For $w\in W$ let $<w>$ be the conjugacy class of $w$.

The elements of $\Irr(W)$ are listed as follows.

(c) $\{1_0,2_1,2_2,\do,2_{(p-1)/2},1_p\}$ ($p$ odd),

(d) $\{1_0,2_1,2_2,\do,2_{(p-2)/2},1',1'',1_p\}$ ($p$ even).

The notation $d_b$ is similar to that in 4.3.
Let $\z=e^{2\pi i/p}$. We have

$A_{<1>,E}=(v^2-1)^2$ if $\dim(E)=1$,

$A_{<1>,E}=2(v^2-1)^2$ if  $\dim(E)=2$,

$A_{<s>,E}= v^2(v^2-1)$ if $E=1_0$,

$A_{<s>,E}= (v^2-1)^2$ if $\dim(E)=2$,

$A_{<s>,E}=-(v^2-1)$ if $E=1_p$,

$A_{<s>,E}=v^2(v^2-1)$ if $E=1'$ ($p$ even),

$A_{<s>,E}=-(v^2-1)$ if $E=1''$ ($p$ even),

$A_{<t>,E}= v^2(v^2-1)$ if $E=1_0$ ($p$ even),

$A_{<t>,E}= (v^2-1)^2$ if $\dim(E)=2$ ($p$ even),

$A_{<t>,E}=-(v^2-1)$ if $E=1_p$ ($p$ even),

$A_{<t>,E}=-(v^2-1)$ if $E=1'$ ($p$ even),

$A_{<t>,E}=v^2(v^2-1)$ if $E=1''$ ($p$ even),

$A_{<c^j>,E}= v^{4j}$ if $E=1_0$,

$A_{<c^j>,E}=(\z^{jk}+\z^{-jk})/p$ if $E=2_k$,

$A_{<c^j>,E}=1$ if $E=1_p$,

$A_{<c^j>,E}=(-1)^jv^{2j}$ if $E=1'$ or $E=1''$ ($p$ even).

\subhead 5.2 \endsubhead
We have

$A'_{1_0,1_0}=1$,

$A'_{2_j,2_k}=(2-\z^{jk}-\z^{-jk})/p$,

$A'_{1_p,1_p}=1$,

$A'_{2_j,1'}=A'_{2_j,1''}=A'_{1',2_j}=A'_{1'',2_j}=(1-(-1)^j)/p$
($p$ even),

$A'_{1',1'}=A'_{1'',1''}=(1-(-1)^{p/2}+p)/2p$ ($p$ even),

$A'_{1',1''}=A'_{1',1''}=(1-(-1)^{p/2}-p)/2p$ ($p$ even).

For all other $E,E'$ we have $A'_{E,E'}=0$.
(The values of $A'_{E,E'}$ are taken from \cite{L94} or, when $p=5$,
from \cite{BM,7.3}.)

\subhead 5.3\endsubhead
Let $h=(v^{2p}-1)(v^4-1)$. We have

$A''_{1_0,1_0}=A''_{1_p,1_p}=v^{2p}h\i$,

$A''_{1_0,1_p}=A''_{1_p,1_0}=h\i$,

$A''_{1_0,2_j}=A''_{2_j,1_0}=A''_{1_p,2_j}=A''_{2_j,1_p}
=(v^{2j}+v^{2p-2j})h\i$,

$A''_{2_j,2_k}=(v^{2j+2k}+v^{2p-2j-2k}+v^{2|j-k|}+v^{2p-2|j-k|})h\i$,

$A''_{1_0,1'}=A''_{1_0,1''}=A''_{1',1_0}=A''_{1'',1_0}=v^ph\i$
($p$ even),

$A''_{1_p,1'}=A''_{1_p,1''}=A''_{1',1_p}=A''_{1'',1_p}=v^ph\i$,
($p$ even),

$A''_{2_j,1'}=A''_{2_j,1''}=A''_{1',2_j}=A''_{1'',2_j}
=(v^{2j+p}+v^{p-2j})h\i$ ($p$ even),

$A''_{1',1'}=A''_{1'',1''}=v^{2p}h\i$ ($p$ even),

$A''_{1',1''}=A''_{1'',1'}=h\i$ ($p$ even).

We shall write $[i,j,k,l]=v^i+v^j+v^k+v^l$.

For $p=7$ the matrix $(A''_{E,E'})$ is
$$h\i\left(\mat
v^{14}& v^2+v^{12}& v^4+v{10}& v^6+v^8& 1\\
v^2+v^{12}&[4,10,0,14]&[6,8,2,12]&[6,8,4,10]&
v^2+v^{12}\\  
v^4+v^{10}&[6,8,2,12]&[6,8,0,14]&[4,10,2,12]&
v^4+v^{10}\\
v^6+v^8&[6,8,4,10]&[4,10,2,12]&[2,12,0,14]&
v^6+v^8 \\
1& v^2+v^{12}&v^4+v{10}&v^6+v^8&v^{14}\endmat\right)$$

For $p=8$ the matrix $(A''_{E,E'})$ is
$$h\i\left(\mat
v^{16}&v^2+v^{14}& v^4+v^{12}& v^6+v^{10}&v^8&v^8&1\\
v^2+v^{14}&[4,12,0,16]&[6,10,2,14]&[8,8,4,12]
&v^6+v^{10}&   v^6+v^{10}&v^2+v^{14}\\
v^4+v^{12}&[6,10,2,14]&[8,8,0,16]&[6,10,2,14]
&v^4+v^{12}&  v^4+v^{12}&v^4+v^{12}\\
v^6+v^{10}&[8,8,4,12]&[6,10,2,14]&[4,12,0,16]
&v^2+v^{14}&  v^2+v^{14}&v^6+v^{10}\\
 v^8&v^6+v^{10}&v^4+v^{12}&v^2+v^{14}&v^{16}&1&v^8\\
 v^8&v^6+v^{10}&v^4+v^{12}&v^2+v^{14}&1&v^{16}&v^8\\
 1&v^2+v^{14}&v^4+v^{12} v^6+v^{10}&v^8&v^8& v^{16}
\endmat \right)$$

\subhead 5.4\endsubhead
The product matrix $AA'$ is as follows.

$(AA')_{<1>,1_0}=(AA')_{1,1_p}=(v^2-1)^2$,

$(AA')_{<1>,2_i}=2(v^2-1)^2$,

$(AA')_{<s>,1_0}=v^2(v^2-1)$,

$(AA')_{<s>,1_p}=-(v^2-1)$,

$(AA')_{<s>,2_i}=(v^2-1)^2$,

$(AA')_{<t>,1_0}=v^2(v^2-1)$ ($p$ even),

$(AA')_{<t>,1_p}=-(v^2-1)$ ($p$ even),

$(AA')_{<s>,1'}=(AA')_{t,1''}=(v^2-1)v^2$ ($p$ even),

$(AA')_{<s>,1''}=(AA')_{t,1'}=-(v^2-1)$ ($p$ even),
       
$(AA')_{<c^j>,1_0}=v^{4j}$,             

$(AA')_{<c^j>,1_p}=1$,

$(AA')_{<c^j>,2_k}=-\d_{jk}v^{2j}$ if $j<p/2$,

$(AA')_{<c^j>,1'}=(AA')_{c^j,1''}=0$ if $j<p/2$ ($p$ even),

$(AA')_{<c^{p/2}>,2_j}=0$    ($p$ even),

$(AA')_{<c^{p/2}>,1'}=(AA')_{<c^{p/2}>,1''}=-v^p$ ($p$ even).

For $p=7$ the matrix $AA'$ is
$$\left(\mat
v^4  &        -v^2&            0&          0&           1\\
v^8  &        0 &           -v^4&          0&           1\\
v^{12} &        0 &           0  &          -v^6&         1\\
v^2(v^2-1)& (v^2-1)^2 &   (v^2-1)^2&    (v^2-1)^2&    -(v^2-1)\\
(v^2-1)^2 &  2(v^2-1)^2& 2(v^2-1)^2& 2(v^2-1)^2&   (v^2-1)^2
\endmat\right)$$
For $p=8$ the matrix $AA'$ is

  $$\left(\mat
v^4 &      -v^2&  0 &  0 &  0 &  0&    v^0\\
v^8 &      0 &   -v^4& 0 &  0 &  0&    v^0\\
v^{12} &     0 &    0&  -v^6&  0&   0&   v^0\\
v^{16}&      0&     0&    0& -v^8& -v^8& v^0\\
  v^2(v^2-1)&  (v^2-1)^2& (v^2-1)^2& (v^2-1)^2&(v^2-1)v^2& -(v^2-1)&     -(v^2-1)\\
  v^2(v^2-1)& (v^2-1)^2& (v^2-1)^2& (v^2-1)^2&  -(v^2-1)& v^2(v^2-1)&     -(v^2-1)\\
  (v^2-1)^2& 2(v^2-1)^2& 2(v^2-1)^2& 2(v^2-1)^2& (v^2-1)^2&(v^2-1)^2&
  (v^2-1)^2 \endmat \right)$$
  
\subhead 5.5\endsubhead
Recall that $\Ps=AA'A''$ is the product of the matrices $AA',A''$.
For $C\in cl(W)$ we write $\Ps_C=\sum_{E\in\Irr(W)}\Ps_{C,E}E$
(a formal sum). We have

$\Ps_{<c>}=1_0$,

$$\align&\Ps_{<c^j>}=(v^{4j-4}+v^{4j-8}+\do+v^4+1)1_0\\&+
\sum_{k;1\le k<j}(v^{4j-4-2k}+v^{4j-8-2k}+\do+v^{2k})2_k\endalign$$
if $1<j<p/2$,
$$\align&\Ps_{<c^{p/2}>}=(v^{2p-4}+v^{2p-8}+\do+v^4+1)1_0\\&+
\sum_{k;1\le k<p/2}(v^{2p-4-2k}+v^{2p-8-2k}+\do+v^{2k})2_k\endalign$$
($p$ even),

$\Ps_{<s>}=\sum_{E\in\Irr(W);E\ne 1_p}E$ ($p$ odd),

$\Ps_{<s>}=\sum_{E\in\Irr(W);E\ne 1_p,E\ne 1''}E$ ($p$ even),

$\Ps_{<t>}=\sum_{E\in\Irr(W);E\ne 1_p,E\ne 1'}E$ ($p$ even),

$\Ps_{<1>}=\sum_{E\in\Irr(W)}\dim(E)E$.

For $p=5$ the matrix $\Ps_{C,E}$ is
$$\left(\mat
1&0&0&0\\
v^4+1&v^2&0&0\\
1&1&1&0\\
1&2&2&1\endmat\right)$$
   
For $p=7$ the matrix $\Ps_{C,E}$ is
$$\left(\mat
1&0&0&0&0\\                                                                v^4+1&v^2&0&0&0\\                                                      v^8+v^4+1&v^6+v^2&v^4&0&0\\
                 1&1&1&1&0\\
                 1&2&2&2&1 \endmat\right)$$
For $p=8$ the matrix $\Ps_{C,E}$ is
$$\left(\mat
1&0&0&0&0&0\\                                                                v^4+1&v^2&0&0&0&0&0\\
v^8+v^4+1&v^6+v^2&v^4&0&0&0&0\\
v^12+v^8+v^4+1&v^10+v^6+v^2&v^8+v^4&v^6&0&0&0\\
1&1&1&1&1&0&0\\
1&1&1&1&0&1&0\\
1&2&2&2&1&1&1 \endmat \right)$$

\subhead 5.6\endsubhead
From the results in 5.5 we see that for our $W$
Theorem 4.2 holds and that
$\Irr_*(W)=\emp,\Irr_{**}(W)=\emp$; moreover,
$\s:cl(W)@>>>\Irr(W)$ is a bijection. It carries the classes of the
elements 5.1(a) (resp. 5.1(b)) in the order written
to the representations in 5.1(c) (resp. 5.1(d)) in the order written,
assuming that $p$ is odd (resp. $p$ is even).

\head 6. Examples\endhead
\subhead 6.1\endsubhead
In this subsection we assume that $W$ is of type $B_3$.
The matrix $\Ps$ is
$$\left(\mat
 1&0&0&0&0&0&0&0&0&0\\
 v^4+1&v^2&0&0&0&0&0&0&0&0\\           
 1&1&1&0&0&0&0&0&0&0\\                 
  1&1&0&1&1&0&0&0&0&0\\                
 1&1&1&1&0&1&0&0&0&0\\                 
A&B&v^8+v^4&v^8+v^4&-2v^6&v^6&0&0&0&0
\\                                                   
v^4+1&v^4+v^2+1&v^4+1&v^2&0&v^2&v^2&0&0&0\\  
 1&2&1&2&1&1&1&1&0&0\\                       
  1&2&2&1&0&2&1&0&1&0\\                       
 1&3&2&3&1&3&3&2&1&1                          
\endmat\right)$$
where
$$A=v^{12}+v^8+v^4+1,B=v^{10}+v^6+v^2.$$

The rows correspond to the conjugacy classes
$$.3,.21,1.2,3.,2.1,.111,1.11,21.1,11.1,111$$ (in this order; notation
of \cite{M15}).
The columns correspond to the irreducible representations
$$1_0,3_1,2_2,3_2,1_3,3_3,3_4,2_5,1_6,1_9$$ (in this order; notation
similar to that in 4.3).
In this case, $\Irr_*(W)$ consists of a single object: $1_3$.

\subhead 6.2\endsubhead
In this subsection we assume that $W$ is of type $H_3$.
The matrix $\Ps$ is
  $$\left(\mat
1&0&0&0&0&0&0&0&0&0\\           
v^4+1&v^2&0&0&0&0&0&0&0&0\\     
1&1&1&1&0&0&0&0&0&0\\           
U&V&v^8+v^4&0&v^6&0&0&0&0&0\\
1&1&1&1&1&1&0&0&0&0\\                           
1&1&2&1&1&1&1&0&0&0\\                   
A&B&C&D&E&F&v^{14}+v^{10}&-2v^{12}&0&0\\  
v^4+1&v^4+v^2+1&v^4+v^2+1&v^4+1&v^2&v^2&v^2&v^2&0&0\\    
1&2&3&2&2&2&2&1&1&0\\                                   
1&3&5&3&4&4&5&3&3&1                           
\endmat\right)$$
where
$$U=v^{12}+v^8+v^4+1,V=v^{10}+v^6+v^2,$$
$$A=v^{24}+v^{20}+v^{16}-v^{12}+v^8+v^4+1,$$
$$B=v^{22}+v^{18}-v^{14}-v^{10}+v^6+v^2,$$ 
$$C=v^{20}+2v^{16}+2v^{12}+2v^8+v^4,$$
$$D=v^{18}-v^{14}-v^{10}+v^6,$$ 
$$E=v^{18}+v^{14}+v^{10}+v^6,F=v^{16}-v^{12}+v^8.$$  
The rows correspond to the conjugacy classes
$$c_3,c_5,(12),c_9,(23),(13),c_{15},(1212),(1),(-)$$
(in this order; notation of 4.4).
The columns correspond to the irreducible representations
$$1_0,3_1,5_2,3_3,4_3,4_4,5_5,3_6,3_8,1_{15}$$ (in this order).

\head 7. Cross sections\endhead
\subhead 7.1\endsubhead
Let $W$ be as in 0.1. We write $W=W_1\T W_2\T\do\T W_n$
where $W_1,W_2,\do,W_n$ are irreducible Coxeter groups.
We can identify $cl(W)=cl(W_1)\T cl(W_2)\T\do\T cl(W_n)$,
$\Irr(W)=\Irr(W_1)\T \Irr(W_2)\T\do\T \Irr(W_n)$ in an obvious
way. We define ${}'\Ph:cl(W)@>>>\Irr(W)$ by
$$(C_1,C_2,\do,C_n)\m({}'\Ph_1(C_1),{}'\Ph_2(C_2),\do,
{}'\Ph_n(C_n))$$
where ${}'\Ph_i(C_i)={}'\Ph(C_i)$ if $W_i$ is a Weyl group and
${}'\Ph_i(C_i)=\s(C_i)$ if $W_i$ is a noncrystallographic
Coxeter group. (Here $C_i\in cl(W_i)$.)

For $C\in cl(W)$ we set
$m(C)=m(w)$ where $w\in C$ and $|C|=|w|$ where $w\in C_{min}$.
For $E_i\in{}'\Ph_i(cl(W_i))$ we set
$${}'X(E_i)=
\{C'\in{}'\Ph_i\i(E_i);m(C')=\min(m(C);C\in{}'\Ph_i\i(E_i)),$$
$${}''X(E_i)=
\{C''\in{}'\Ph_i\i(E_i);m(C'')=\max(m(C);C\in{}'\Ph_i\i(E_i)).$$
We define
$\t':{}'\Ph(cl(W))@>>>cl(W)$, $\t'':{}'\Ph(cl(W))@>>>cl(W)$ by
$$\t'(E_1,E_2,\do,E_n)=(C'_1,C'_2,\do,C'_n),$$
$$\t''(E_1,E_2,\do,E_n)=(C''_1,C''_2,\do,C''_n)$$
(with $E_i\in\Irr(W_i),C'_i\in{}'X(E_i),C''_i\in {}''X(E_i)$)
where

(a) $|C'_i|=\max(|C|;C\in {}'X(E_i))$,

(b) $|C''_i|=\min(|C|;C\in {}''X(E_i))$.
If $W_i$ is a Weyl group, each of the sets ${}'X(E_i),{}''X(E_i)$
consists of one element (see \cite{L15}, \cite{L22}) so that
(a),(b) are automatically satisfied. If $W_i$ is
noncrystallographic, the sets ${}'X(E_i),{}''X(E_i)$ may have more
than one element; but imposing the conditions (a),(b) guarantees
that $C'_i$ and $C''_i$ are uniquely determined. Thus the maps
$\t',\t''$ are well defined. We have

${}'\Ph\t'(E)=E$, ${}'\Ph\t''(E)=E$ for all $E\in{}'\Ph(cl(W))$.

\subhead 7.2\endsubhead
If $W$ is of type $H_4$, the image of $\t'$ has $22$ elements:
$$\align&c_4,c_6,c_8,c_{10},c_{12},c_{16},c'_{16},c_{18},c_{20},c_{22},
\\&c_{24},c_{26},(124),c_{28},c_{30},(12),c_{40},(23),c_{48},c_{60},
(1),(-);\endalign$$
the image of $\t''$ has $22$ elements:
$$\align&c_4,c_6,c_8,c_{10},c_{12},(123),c'_{16},c_{18},c_{20},
(12123),c_{24},\\&c_{26},(124),c_{28},(243),(12),
(134),(23),(12124),(13),(1),(-).\endalign$$
                     
If $W$ is of type $H_3$, the image of $\t'$ has $8$ elements:

$c_3,c_5,(12),c_9,(23),c_{15},(1),(-);$

the image of $\t'$ has $8$ elements:

$c_3,c_5,(12),c_9,(23),(1212),(1),(-)$.

If $W$ is a dihedral group, $\t'=\t''$ is the bijection
inverse to ${}'\Ph$.

\enddocument